\documentclass[11pt]{article}
\usepackage[a4paper,margin=27mm]{geometry}
\usepackage{amsmath,amssymb,amsthm,hyperref,xcolor}
\newcommand{\added}[1]{{\color{black}#1}}
\newcommand{\reviewfix}[1]{{\color{black}#1}}
\hypersetup{colorlinks=true,linkcolor=black,urlcolor=black,citecolor=black}
\newcommand{\R}{\mathbb R}
\newcommand{\BMO}{\operatorname{BMO}}
\newcommand{\avint}{\mathop{\mathrm{Avg}}}
\newcommand{\loc}{\mathrm{loc}}
\newcommand{\dist}{\operatorname{dist}}
\newtheorem{theorem}{Theorem}
\newtheorem{lemma}{Lemma}

\title{BMO preservation and quasiconformality
}
\author{Yepei Wen and Yuan Zhou\thanks{Corresponding author: Yuan Zhou.}
}

\begin{document}
	\maketitle
	\begin{abstract}  
		In dimension  $n\ge2$, a homeomorphism of $\R^n$  preserving BMO space therein was proved to be   quasiconformal\reviewfix{, up to an orientation reversing reflection}. 
	\end{abstract}
	\noindent\textit{Mathematics Subject Classification (2020):}
	Primary 30C65; Secondary 42B35, 46E35.\par
	\vspace{-1.2em}
	\section{Introduction}
	Let $n\ge 2$ and denote by  BMO the space of functions having bounded mean oscillations on $\mathbb R^n$, that  is,   the collection of all $u\in L^1_\loc(\mathbb R^n)$ with the seminorm  
	\[
	[u]_{\BMO}=\sup_B\avint_B|u-u_B|<\infty,
	\]
	where the supremum is taken over Euclidean balls. Throughout this paper,  write 
	$$\avint_Eu=u_E=|E|^{-1}\int_Eu(z)\,dz $$
	\added{for every bounded measurable set $E$ with $|E|>0$.}  Observe that  translation scaling and rotation induce  isometries of BMO. Indeed, if $f$ is one of 
	translation scaling and rotation, then $[u\circ f]_{\BMO}=[u]_{\BMO} $ for all $u\in \BMO$.

	In 1974\added{,} Reimann  \cite{Reimann} considered \added{the} \added{quasiconformal} invariance of BMO. Recall that a \added{homeomorphism} $f:\mathbb R^n\to\R^n$ is called $K$-quasiconformal for some $K\ge 1$ if $f\in \dot W^{1,n}_\loc(\mathbb R^n)$ and 
	$$|Df|^n\le K\det Df\ \text{\added{almost everywhere in}}\ \mathbb R^n.$$
	\added{Here $|Df|$ is the operator norm. Composition in \eqref{y1} is understood for every finite measurable representative; on almost-everywhere classes it must be well defined.}
	Reimann  \cite{Reimann} proved that 
	\begin{align}\label{y1}
		[u\circ f]_{\BMO}\le C[u]_{\BMO}\quad\forall u\in \BMO.
	\end{align}
	for some constant $C=C(n,K)$.
	
	Since then, the following converse problem is left open in the literature; see for example \cite{KMS}. 
	
	\noindent {\bf Question.}  If $f$ is  a homeomorphism of $\R^n$ satisfying   \eqref{y1}, then  is $f$ a $K$-\added{quasiconformal} mapping for some constant  $K\ge 1$, allowing an orientation reversing reflection\reviewfix{?}

	Towards this,  if  $f$ is almost everywhere differentiable and 
	\added{an ACL regular homeomorphism} of $\mathbb R^n$, and induces \added{a} BMO isomorphism, 
	Reimann 1974 proved that $f$ is a quasiconformal mapping, allowing an orientation reversing reflection.   Assuming  a homeomorphism $f$ of $\R^n$ 
	\added{preserves} all BMO spaces in   subdomains of $\mathbb R^n$, Astala \added{\cite{Astala}} proved that $f$ must be $K$-quasiconformal, allowing an orientation reversing reflection. 
	Gotoh \cite{Gotoh}  obtained relative density property for homeomorphisms of $\R^n$ satisfying  \eqref{y1}; see also  \cite{KKMS}. 
	
	The aim of this paper  is to give a full positive answer for the above question. 	
	\begin{theorem}\label{main}
		Let $n\ge2$. If $f:\R^n\to\R^n$ \added{is} a homeomorphism 
		satisfying \eqref{y1}, then    $f$ is $K$-quasiconformal for some constant  $K$ depending on \added{$n$ and} the constant $C$ therein, allowing an orientation reversing reflection.		 
	\end{theorem}
	
	The proof starts from the
	relative density property in Lemma~\ref{2} established by
	Gotoh~\cite{Gotoh}.
	It  implies a  quantitative upper bound for
	distance of two sets $E_1,E_2$ contained in \reviewfix{$f(B)$} via  $\mu(B)^{1/n}$  and lower bound of densities of their pre-images; see  Lemma 3.  Here and below, define $\mu(E)=|f(E)|$ for Borel sets $E$.
	For each ball $B$, Lemma 3 allows to choose vectors $m_B$   from certain \reviewfix{coordinate quantiles} of $f$ in $B$, and find large set $G_B\subset B$ so that  $|m_B-f(z)|$ is bounded via $\mu(B)^{1/n}$; see Lemma 4. 
	
	Since $m_{B(x,r)}$ converges to $f(x)$, Lemma 4 and a telescope argument give  the pointwise estimate of Haj\l asz type 
	$$|f(x)-f(y)|\le C\reviewfix{|x-y|}\{[M_R(\mu)(x)]^{1/n}+[M_R(\mu)(y)]^{1/n}\}$$
	whenever \reviewfix{$x,y\in B(0,R)$ with $0<|x-y|<R$}; see Lemma~\ref{pointwise}.
	Here $M_R(\mu)$ is the local Hardy maximal function
	of the measure $\mu$. 
	Applying the weak $(1,1)$ estimate,\reviewfix{$[M_R(\mu)]^{1/n}$}$\in 
	L^{ p}_{\loc}(\R^n)$ for $1<p<n$.         
	The Haj\l asz--Sobolev characterization then gives
	$f\in W^{1,p}_{\loc}(\R^n)$ for $1<p<n$; see Lemma~\ref{subcrit}.
	Moreover, applying the  quantitative upper bound for distances of suitable two subsets arising from approximate  differentiability, one further shows that 
	$|Df|^n$ is bounded \reviewfix{by a constant times the Radon--Nikodym density of the absolutely continuous part of $\mu$} and hence  $f\in W^{1,n}_{\loc}(\R^n)$; see Lemma~\ref{critical}.

	Finally, we adapt Reimann's argument using approximate differentiability.
	Lemma~\ref{rescaling} gives local $L^1$ convergence of bounded
	continuous test functions composed with the rescaled maps to those
	composed with $z\mapsto Df(x)z$.
	This transfers the BMO composition bound to $z\mapsto Df(x)z$
	at almost every $x$, for bounded continuous BMO test functions.
	The logarithmic test functions of Lemma~\ref{tests} then show that
	for almost every $x$, either $Df(x)=0$ or $Df(x)$ is invertible
	with the ratio of its largest to smallest singular value bounded
	in terms of $n$ and the constant $C $ in (1); see Lemma~\ref{affine}.
	Consequently,
	\[
	|Df|^n\le C|\det Df|\qquad\text{almost everywhere}.
	\]
	Together with $W^{1,n}_{\loc}$ regularity and the almost-everywhere
	sign of the Jacobian determined by the topological orientation,
	this proves quasiconformality, up to a fixed orientation reversing
	reflection.

	{\bf Acknowledgment.} Since 2013, Yi Ru-Ya Zhang and Yuan Zhou started to work on this question. At that moment, 
	they realized   that  differentiability used in Reimann's proof could be reduced to \added{approximate differentiability} so to get: 
	If a Sobolev $\dot W^{1,1}_\loc$-\added{homeomorphism} of $\R^n$ induces an isomorphism of BMO, then it must be quasiconformal quantitatively.  This was unpublished. 
	Essentially, Lemma \ref{rescaling}-\ref{affine} below is due to Yi Ru-Ya Zhang and Yuan Zhou. 
	
	Recently, Yepei Wen and Yuan Zhou asked \added{OpenAI} the possibility to get the 
	\added{Sobolev regularity of homeomorphisms}  preserving  BMO space. 
	\added{OpenAI} gave the key idea to derive sobolev regularity from Lemma \ref{2}. 
	Yepei Wen and Yuan Zhou  follow this idea  and 
	give all necessary details so to get Lemma \ref{sep}-\ref{critical} below.

	Yepei Wen and Yuan Zhou are responsible for the whole proofs in this paper.
	They would like to thank Yi Ru-Ya Zhang for his important contributions to this paper mentioned as above. 

	\section{Proof of Theorem \ref{main}}
	
	Assume that $f:\R^n\to\R^n$ is a homeomorphism satisfying 
	$$\|u\circ f\|_{\BMO(\R^n)}\le C_f\|u\|_{\BMO(\R^n)}\quad\forall u\in\BMO(\R^n)$$
	\added{where} the constant $C_f\ge1$. \added{For $j=1,\cdots,n$, denote by $f_j$   the $j$th coordinate of $f$.} 
	Define $\mu(E)=|f(E)|$ for all Borel sets $E\subset\R^n$. \added{This is a locally finite Radon measure. Write $\omega_n=|B(0,1)|$.}  
	\begin{lemma}\label{bounded} 
		The inverse map  $f^{-1}$ satisfies the Lusin N-property, that is, 
		$|E|=0$ implies $|f^{-1}(E)|=0$. Consequently, 		 
		for each $j$, \added{a Borel null set $E\subset\R$ satisfies $|f_j^{-1}(E)|=0$.}
	\end{lemma}
	\begin{proof}  
		Assume that  $f^{-1}$ does not satisfy the Lusin N-property. Then there is a   Borel null set $E$
		so that $|f^{-1}(E)|>0$.   Choose a
		compact set $K\subset f^{-1}(E)$ with $0<|K|<\infty$. 
		Write $K$ as the union of a family $\{K_j\}$ of mutually disjoint Borel sets with positive measures.
		Define  	$$h=\sum_{j=1}^\infty\frac{\chi_{K_j}}{j|K_j|} .$$
		Obviously, $$\int_Kh\,dx=\sum_{j=1}^\infty \frac1j=\infty$$ 
		Since  $|f(K)|\le |E|=0$, the function  
		$u=( h\circ f^{-1})\chi_{f(K)}$ is measurable and zero almost everywhere, 
		so  $u$  is the $0$-function in BMO.
		But $u\circ f=h$ on $K$ is not locally integrable, and hence not in BMO. 
		This is a contradiction.

		Finally, note that $ f_j^{-1}(E)  =f^{-1} (\{x\in\R^n|x_j\in E \})$ for any Borel set $E\subset\R$. 
		If $|E|=0$, then  $|\{x\in\R^n|x_j\in E \}|=0$ and hence $| f_j^{-1}(E)|=0$.
		
	\end{proof}
	The following was proven by \added{Gotoh \cite{Gotoh}}.  
	For Borel sets $E_1,E_2\subset\R^n$, define 
	\[
	\Delta(E_1,E_2)=\sup_B\min_{i=1,2}\frac{|E_i\cap B|}{|B|}.
	\]
	\begin{lemma}\label{2}
		There  are constants
		$C,\alpha>0$ depending on $n,C_f$ such that, for every pair $E_1,E_2$
		of Borel sets in $\R^n$,
		\begin{equation}\label{x0}
			\Delta(f^{-1}(E_1),f^{-1}(E_2))
			\le C\Delta(E_1,E_2)^\alpha.  
		\end{equation} 
	\end{lemma}  
	
	The following is a direct  consequence of Lemma \ref{2}. 
	\begin{lemma} \label{sep}
		Let $\delta\in(0,1/2)$.   Given any ball $B$, and Borel sets $E_1,E_2\subset f(B)$,
		if 
		\begin{equation}\label{x1}
			|B\cap f^{-1}(E_i)|\ge\delta|B| \qquad \text {for}\ i=1,2,
		\end{equation}
		then
		\begin{equation}\label{x2}
			\dist(E_1,E_2)\le \added{C(n,C_f,\delta)}\mu(B)^{1/n}. 
		\end{equation}
	\end{lemma}
	
	\begin{proof} Assume that $\dist(E_1,E_2)>0$ since the case  $\dist(E_1,E_2)=0$ is trivial.  Note that 
		\begin{align*}
			\Delta(E_1,E_2)  
			&=
			\sup\left\{
			\min_{i=1,2}\frac{|E_i\cap B(z,r)|}{|B(z,r)|}|  B(z,r)\cap E_i\ne \emptyset,i=1,2 \right \}
		\end{align*}
		Given any ball $B(z,r)$, which has nonempty intersection with both $ E_1$ and $E_2$,
		observing   $r>\frac12\dist(E_1,E_2)$,  
		by $E_1,E_2\subset f(B)$  we obtain 
		\begin{align*}
			\Delta(E_1,E_2)   \le \frac{2^n}{\omega_n}\frac{|f(B)| }{\dist(E_1,E_2)^n} \end{align*}
		Since the assumption \eqref{x1} gives $$\Delta(f^{-1}(E_1),f^{-1}(E_2))\ge\delta,$$ 
		applying \eqref{x0} we obtain  the desired \eqref{x2}.
	\end{proof}
	
	From Lemma \ref{sep}, we deduce the following. 
	\begin{lemma} \label{quant}  
		Fix $0<\varepsilon<1/2$.
		For every ball $B\subset\R^n$, there are a vector $m_B\in\R^n$ and a Borel set $G_B\subset B$
		such that
		\[
		|B\setminus G_B|\le\varepsilon|B| \qquad\text{and}\ 
		\sup_{z\in G_B}|f(z)-m_B|\le \added{C(n,C_f,\varepsilon)}|\mu (B)|^{1/n} .
		\tag{Q}
		\]
		Moreover, $m_{B(x,r)}\to f(x)$ as $r\downarrow0$ at all $x\in\R^n$.
	\end{lemma}
	\begin{proof}      	
		\added{Lemma~\ref{bounded} gives $|B\cap f_j^{-1}(\{t\})|=0$ for all $t$, so the distribution function of $f_j$ on $B$ is continuous.} For each $j$, choose $a_j\le b_j$ such that 
		\[
		|\{z\in B:f_j(z)<a_j\}|=\frac\varepsilon {2n}|B|,\qquad
		|\{z\in B:f_j(z)>b_j\}|=\frac\varepsilon {2n}|B|.
		\]
		Set    $$ G_B=\added{B\cap}\bigcap_{j=1}^nf_j^{-1}([a_j,b_j]).$$          Obviously,  $$|B\setminus G_B |=|\bigcup _{j=1}^n B\setminus f_j^{-1}([a_j,b_j])|\le \added{\varepsilon}|B|. $$

		Next, write  \added{$E_1=f(B)\cap\{y:y_j<a_j\}$ and $E_2=f(B)\cap\{y:y_j>b_j\}$}. Then $\dist(E_1,E_2)\ge |b_j-a_j|$.
		Lemma \ref{sep} gives
		\[
		b_j-a_j\le \dist(E_1,E_2)\le C(\added{\varepsilon}) |f(B)|^{1/n}.
		\] 
		Moreover, set   $$m_B=\left(\frac{b_j+a_j}2\right)_{j=1}^n.$$
		For any $z\in G_B$, we have $f_j(z)\in [a_j,b_j]$ for all $j$, 
		and hence $$|f(z)-m_B|=\left(\sum_{j=1}^n|f_j(z)- \frac{b_j+a_j}2|^2\right)^{1/2}\le \added{\frac12\big(\sum_{j=1}^n|b_j-a_j|^2\big)^{1/2}} \le C(\added{\varepsilon}) |f(B)|^{1/n}$$
		as desired. 
		
		Finally, since $$\inf_B f_j\le (m_B)_j\le \sup_B f_j\ \text{ for all}\  j,$$ 
		the continuity of $f$ gives $m_{B(x,r)}\to f(x)$ as $r\to0$ for all $x\in\R^n$.
		
	\end{proof}

	Lemma \ref{quant} allows to get the following pointwise estimate of Hajlasz type. 
	The local Hardy-Littlewood maximal function of $\mu$ is given by 
	$$M_R{\mu }(x)=\sup_{x\in B(z,r),r<R}\frac{ \mu (B(z,r))}{|B(z,r)|}.$$

	\begin{lemma} \label{pointwise}
		\reviewfix{There is a constant} $C(n,C_f)$   such that,
		\[
		|f(x)-f(y)|\le C|x-y|\big[ \big[M_R(\mu )(x)\big]^{1/n} + \big[M_R (\mu )(y)\big]^{1/n} \big]. \]
		for all  $x,y\in B(0,R)$ with $0<|x-y|< R $. 
	\end{lemma}
	\begin{proof} 
		Fix any $x, y \in B(0,R)$ with $ r=|x-y|< R$.  
		We claim that 
		\begin{align}\label{claim-0}
			|f(x)-m_{B(x,r)}| \le C r \big[M_R(\mu )(x)\big]^{1/n}.
		\end{align}
		To see this, set $B_k=B(x,2^{-k}r)$ for $k\ge 0$. By Lemma \ref{quant},
		$$
		|f(x)-m_{B(x,r)}|=\added{\lim_{k\to\infty}|m_{B_k}-m_{B_0}|\le\sum_{k=0}^\infty} |m_{B_k}-m_{B_{k+1}}|.$$
		It then suffices to prove 
		\begin{align}\label{claim}
			|m_{B_k}-m_{B_{k+1}}|\le C 2^{-k}r \added{\big(M_R(\mu)(x)\big)^{1/n}}.
		\end{align}
		For each $k$, denote by  $G_{B_k}$ the set given in Lemma \ref{quant} with $\varepsilon=2^{-n-4}$. 
		Since $B_{k+1}\subset B_k$, by Lemma \ref{quant} we have 
		\begin{align*}
			|B_{k+1} \setminus[G_{B_k}\cap G_{B_{k+1}}]| &=  | (B_{k+1} \setminus G_{B_k})\cup 
			(B_{k+1} \setminus G_{B_{k+1}})|\\
			&\le 
			| (B_k \setminus G_{B_k})\cup 
			(B_{k+1} \setminus G_{B_{k+1}})|\\
			&\le 2^{-n-4}|B_k|+ 2^{-n-4}|B_{k+1}|\\
			&<2^{-3}|B_{k+1}|. 
		\end{align*}
		Thus  		\[
		|G_{B_k}\cap G_{B_{k+1}}|= |B_{k+1}| -|B_{k+1} \setminus[G_{B_k}\cap G_{B_{k+1}}]|
		\ge |B_{k+1}|- 2^{-3}|B_{k+1}|>0.
		\]
		Choosing $z\in G_{B_k}\cap G_{B_{k+1}}$, by Lemma \ref{quant} again we obtain
		\begin{align*}
			|m_{B_k}-m_{B_{k+1}}|&\le |m_{B_k}-f(z)|+|f(z)-m_{B_{k+1}}|\\
			& 		\le C\big[\mu(B_k)^{1/n}+\mu(B_{k+1})^{1/n}\big]\\
			&		\le C \mu(B_k)^{1/n}.
		\end{align*}
		Since 		\[
		\mu(B_k) 
		\le \omega_n 2^{-kn}r^n\,M_R(\mu )(x),
		\]
		we get  \eqref{claim} and hence \eqref{claim-0}. 
		
		Changing the role of $x,y$ we know that 
		$$
		|f(y)-m_{B(y,r)}| \le C r [M_R(\mu )(y)]^{1/n}.$$
		
		Finally,         write $W=[B(x,r)\cap B(y,r)]$.
		Since  
		$B(\frac{x+y}2,\frac r2)\subset W$, we have  $|W|\ge 2^{-n}|B(x,r)|$.
		Lemma \ref{quant} with $\varepsilon=2^{-n-4}$ then gives 
		\begin{align*} |W\setminus [G_{B(x,r)}\added{\cap}G_{B(y,r)} ]|&=
			|[W\setminus G_{B(x,r)}]\added{\cup}[W\setminus
			G_{B(y,r)}]|\\ 
			&\le |B(x,r)\setminus  G_{B(x,r)}|+|B(y,r) \setminus
			G_{B(y,r)}|\\
			& \le 2^{-n-3}|B(x,r)|\\
			&\le \frac14 |W|.\end{align*}
		Thus  $$|W\cap[G_{B(x,r)}\added{\cap}G_{B(y,r)}]|\ge \frac34 |W|>0.$$ 
		Choosing  $ z\in G_{B(x,r)}\cap G_{B(y,r)} $, by Lemma \ref{quant} we obtain 
		\begin{align*}
			|m_{B(x,r)}-m_{B(y, r)}|&\le|m_{B(x,r)}-f(z)|+|f(z)-m_{B(y, r)}|\\
			&\le \added{C\big[\mu(B(x,r))^{1/n}+\mu(B(y,r))^{1/n}\big]} \\
			& \le Cr [ (M_R(\mu)(x))^{1/n}+ (M_R(\mu)(y))^{1/n}]
		\end{align*}
		Since
		$$|f(x)-f(y)|\le |f(x)-m_{B(x,r)}|+  |m_{B(x,r)}-m_{B(y, r)}|+|f(y)-m_{B(y,r)}|,$$
		combining above estimates together we get the desired upper bound. 
	\end{proof}
	
	Applying the weak estimate of maximal function of measures appeared in Lemma \ref{pointwise}, we  \added{obtain the following.}  
	\begin{lemma}\label{subcrit}
		We have 
		$f\in W^{1,p}_\loc(\R^n)$ for $1\le p<n$. 
	\end{lemma}
	\begin{proof}
		Thanks to the weak estimates for $M_R[\mu]$, one has
		\[
		\added{|\{x\in B(0,R):(M_R\mu(x))^{1/n}>t\}|
			\le C\mu(B(0,3R))t^{-n}.}
		\]
		\added{Indeed, every ball of radius less than $R$ containing $x\in B(0,R)$ lies in $B(0,3R)$.}
		Applying this for $t>1$, we get
		\begin{align*}
			\int_{B(0,R)}(M_R\mu)^{p/n}\,dx
			&\added{=p}\int_0^\infty t^{p-1}
			|\{x\in B(0,R):(M_R\mu(x))^{1/n}>t\}|\,dt\\
			&\le\added{|B(0,R)|+Cp\mu(B(0,3R))\int_1^\infty t^{p-1-n}\,dt}<\infty
		\end{align*}
		for $1<p<n$, as desired.
		
		This means that $f$ belongs to the Haj\l asz--Sobolev space $  M^{1,p}_\loc(\R^n)$ with $1<p<n$.
		Since \added{Haj\l asz \cite{Hajlasz}} proves $  M^{1,p}_\loc(\R^n)=    W^{1,p}_\loc(\R^n)$,
		we therefore get $f\in W^{1,p}_\loc(\R^n)$ with $1<p<n$. \added{Here the pointwise bound is applied on small balls outside the null set $\{M_R\mu=\infty\}$; the case $p=1$ follows locally by H\"older\textquotesingle s inequality.}
	\end{proof}
	
	Assume that  $ u\in  W^{1,1}_\loc(\R^n)$. 
	Then $u$ \added{admits} an ACL representative $\widetilde u$, which 
	\added{  is absolutely continuous on almost all lines parallel to the axes.}
	We may write $u$ as $\widetilde u$;  if $u$ is continuous, \added{we may take} $\added{\widetilde u}=u$.
	Thus the partial derivatives $Du=(\partial_{x_i}u)_{1\le i \le n}$ exist  almost everywhere
	and coincide with the weak/distributional derivative of \added{$u$} almost everywhere. 
	\added{By the standard approximate differentiability theorem for Sobolev functions}, we know that $u$ is approximately differentiable almost everywhere 
	and its approximate differential coincides with $Du$ almost
	everywhere. In other words, for almost all $x\in\R^n$  
	$$ap.\lim_{z\to0}\frac{|u(x+z)-u(x)-Du(x)\cdot z|}{|z|}=0$$
	where the approximate limit $ap.\lim_{z\to0 } v(z)=l$  means  that 
	for any $\epsilon>0$, 
	$$\lim_{r\to0 }\frac{\added{|\{z\in B(0,r):|v(z)-l|>\epsilon\}|}}{\added{|B(0,r)|}}=0.$$
	If \added{$ap.\lim_{z\to0}v(z)=v(0)$}, we say that $v$ is approximately continuous at \added{$0$}.

	Thanks to Lemma \ref{subcrit}, and continuity of $f$, we know that each $f_j$ is its ACL representative. 
	Thus   at almost all $x\in\R^n$, the approximate differential of $f$ is given by the partial derivatives $Df(x)$. Moreover, write
	$\mu=w\,dx+\mu_s$ \added{with $\mu_s\perp dx$} for some \added{nonnegative} function $w\in L^1_\loc(\R^n)$.
	\added{By the differentiation theorem for Radon measures,} 
	\begin{equation}\label{metric}
		\lim_{r\downarrow0}\frac{\mu(B(x,r))}{|B(x,r)|
		}=w(x) 
	\end{equation}
	for almost all $x\in\R^n$. 
	Denote by $E_f$ the set of all $x$   at which  the approximate differential of $f$ 
	is given by  $ Df(x)$, and  \eqref{metric} holds. 
	We have $|\R^n\setminus E_f|=0$.

	The following follows from the \added{approximate differentiability} and  Lemma \ref{2}.

	\begin{lemma} \label{critical}
		We have 
		\(|Df |^n\le Cw \) almost everywhere in $\R^n$. 
		Consequently,  $f\in W^{1,n}_{\loc}(\R^n)$.
	\end{lemma}
	
	\begin{proof}
		
		
		\added{Fix $x\in E_f$. The case $Df(x)=0$ is immediate.} Below assume $Df(x)\ne 0$.
		\added{There must exist} some $j$ such that $|Df_j(x)|\ge \frac1n |Df(x)|>0$.  
		For any $\delta\in(0,1/4)$, there exists  
		$c(\delta,n)>0$ such that  
		\[\left|\{z\in B(0,r):\pm z_1>2cr \}\right|  = 2\delta |B(0,r)|.\]
		Write 
		\[A_\pm(r)= 
		\{z\in B(x,r)|\pm \frac{Df_j(x)}{|Df_j(x)|}\cdot (z-x) >2cr \}\]
		We have 
		\[ | A_\pm  | 
		=\left|\{z\in B(0,r): z_1 >2cr \}\right|\added{=2\delta |B(0,r)|}\]
		
		Set  
		$$ F_\pm (r)= \{z\in B(x,r)| \added{\mp}[f_j(z)-f_j(x)-Df_j(x)(z-x)]>  c|z-x||Df_j(x)|\}.$$
		For all sufficiently small $r$, approximate differentiability of $f_j$ at $x$ implies that 
		\[
		|\{z\in B(x,r)\reviewfix{:}\ |f_j(z)-f_j(x)-Df_j(x)(z-x)|>  c|z-x||Df_j(x)|\}|\le \delta |B(x,r)|, 
		\] 
		
		If  $z\in A_\pm(r)$ and $z\notin F_\pm(r)$,  
		we have 
		{\color{black}
			\begin{align*}
				\pm[f_j(z)-f_j(x)]
				&=\pm Df_j(x)(z-x)\pm[f_j(z)-f_j(x)-Df_j(x)(z-x)]\\
				&>2cr|Df_j(x)|-c|z-x||Df_j(x)|\\
				&\ge cr|Df_j(x)|.
			\end{align*}
		}
		Write   $$H_{\pm}= \{z\in B(x,r):\pm[f_j(z)-f_j(x)]> cr\added{|Df_j(x)|}\}$$
		we know that $ |H_\pm|\ge  \delta |B(x,r)|. $ 
		
		Since $f(H_\pm)\subset f(B(x,r))$ and $\dist(f(H_+),f(H_-))\added{\ge}2cr|Df_j(x)|$, applying   Lemma \ref{sep} we obtain 		\[
		(2cr|Df_j(x)|)^n\le \added{C(n,\delta,C_f)}\mu(B(x,r)).
		\]
		\added{Dividing by $r^n$ and passing to the limit}, we get   $|Df_j(x)|^n\le\added{C}w(x)$ as desired. \added{Since $|Df(x)|\le n|Df_j(x)|$ and $w\in L^1_\loc$, the existing weak derivatives belong to $L^n_\loc$, proving the last assertion.} 
		
	\end{proof}
	
	
	The following $L^1_\loc(\R^n)$ convergence follows from \added{approximate differentiability}. 
	For any $x\in\R^n$ and $r>0$, write
	$$f_{x,r}(z)=\frac{f(x+rz)-f(x)}r.$$
	Given any $x\in E_f$, put the linear function 
	$$\ell_x(z)=Df(x)z\forall z\in\R^n$$ 
	
	\begin{lemma} \label{rescaling}
		Given any   $x\in E_f$, for any $v\in C^0(\R^n)\cap L^\infty(\R^n)$,  we have 
		$v\circ f_{x,r}(z)\to v\circ \ell_x$ in  $L^1_\loc(\R^n)$    as $r\to0$.
	\end{lemma}
	\begin{proof}
		Fix $R>0$. Write $T=\|v\|_{L^\infty(\R^n)}$. 
		For any $\epsilon>0$, we have 
		\begin{align*}  \int_{B(0,R)}|v\circ f_{x,r}(z)-v\circ\ell_x(z)|\,dz&=\added{\int_0^{2T}|\{z\in B(0,R):|v\circ f_{x,r}-v\circ\ell_x|>t\}|\,dt} \\
			&\le \epsilon |B(0,R)|+ \added{2T}|\{z\in B(0,R)||v\circ f_{x,r}(z)-v\circ\ell_x(z)|>\epsilon\}|
		\end{align*}
		It then suffices to show that the second term converges to $0$ as $r\to0$.
		Thanks to the uniform continuity of $v$ on
		\added{the compact set $\{q:\dist(q,\ell_x(\overline{B(0,R)}))\le1\}$},
		there exists $0<\eta<\added{\min\{1,\epsilon\}}$
		such that
		\added{for $z\in B(0,R)$, if $|f_{x,r}(z)-\ell_x(z)|\le\eta$, both points lie in this compact set and their $v$-values differ by at most $\epsilon$. Hence}
		\begin{align*} &\{z\in B(0,R)||v\circ f_{x,r}(z)-v\circ\ell_x(z)|>\epsilon\}\\
			&\quad\subset 
			\{z\in B(0,R)||  f_{x,r}(z)- \ell_x(z)|>\eta\}\\
			&\quad\subset \cup_{j=1}^n
			\{z\in B(0,R)||  (f_j)_{x,r}(z)-\added{Df_j(x)z}|>\eta/n\}.
		\end{align*}
		For each $j$ one has 
		\begin{align*}
			&|\{z\in B(0,R)| |(f_j)_{x,r}(z)-Df_j(x)z|>\eta/n   \}|\\
			&\quad=|\{z\in B(0,R)| |f_j(x+rz)-f_j(x)-rDf_j(x)z|>\eta r/n \}|\\
			&\quad=\added{r^{-n}}|\{y\in B(x,rR):|f_j(y)-f_j(x)-Df_j(x)(y-x)|>\eta r/n\}|  \\
			&\quad \le r^{-n}|\{y\in B(x,rR):|f_j(y)-f_j(x)-Df_j(x)(y-x)|>
			(\eta/(nR))|y-x|\}| \\
			&\quad =r^{-n}|B(x,rR)|o(1)\\
			&\quad =|B(0,R)|o(1). 
		\end{align*}
		Thus 
		$$
		|\{z\in B(0,R)||v\circ f_{x,r}(z)-v\circ\ell_x(z)|>\epsilon\}|=\added{o(1)}$$
		\added{Letting $r\downarrow0$ and then $\epsilon\downarrow0$ proves the assertion.} 
	\end{proof}

	\added{The following test function} will be used later. Let  $\zeta (s)=\max\{-1,\min\{s,1\}\}$.
	For \added{$N>0$}, write 
	\[L_N(t)=\min\{\added{\log^+}\frac1{|t|},N\}\quad\text{for }\added{t\ne0},
	\qquad L_N(0)=N.
	\]
	and define  
	$$v_N(y)=\zeta(y_1)L_N(y_n)\quad\forall y\in\R^n.$$

	\begin{lemma}\label{tests}
		The function 	$
		v_N \in C^0(\R^n)\cap L^\infty(\R^n)$
		and $[v_N]_{\BMO(\R^n)}\le C(n)$ uniformly
		for $N>0$.
	\end{lemma}
	\begin{proof}
		
		
		Consider any cube $Q=I_1\times\cdots I_n$ with center $x=(x_1,\cdots,x_n)$. 
		Let $a(I_n)=\avint_{I_n}L_N$.  we have 
		\begin{align*}
			|v_N(z)-\zeta(x_1)a(I_n)|&\added{\le}|v_N(z)- \zeta(z_1)a(I_n)|+ |\zeta(z_1)\added{a(I_n)}-\zeta(x_1)a(I_n)| \\
			&
			\le   |L_N(z_n)-a(I_n)|
			+ |a(I_n)| |\zeta(z_1)-\zeta(x_1)|.
		\end{align*}
		and hence 
		\[
		\avint_Q|v_N-\zeta(x_1)a(I_n)|
		\le \avint_{I_n}|L_N-a(I_n)|
		+ |a(I_n)|\avint_{I_1}|\zeta-\zeta(x_1)|.
		\]
		It is known that \added{$[L_N]_{\BMO(\R)}\le C$ uniformly in $N$, since $\log|t|\in\BMO(\R)$ and scalar $1$-Lipschitz truncation increases its seminorm by at most a factor two.} By $|\zeta|\le1$ and $|\zeta^\prime|\le 1$ we have 
		\[
		\avint_Q|v_N-\zeta(x_1)\added{a(I_n)}|
		\le\added{C}
		+ 2 |a(I_n)|\min\{1,|I_1|\}  
		\]
		Note that 
		$$|a(I_n)|\added{\le}\frac1{|I_n|} \int_{I_n\cap[-1,1]} \log \frac1{|t|}\,dt.$$
		If $ |I_n|\ge 1$, we have 
		$$\added{|a(I_n)|}\le  \frac2{|I_n|} \int_0^1\log\frac1t\,dt\le C. $$
		If $ |I_n|<1$, we have 
		$$|a(I_n)| \le  \frac2{|I_n|} 
		\int_0^{|I_n|/2}\log\frac1t\,dt\le\added{1+}\log\frac2{|I_n|}. $$
		By $|I_1|=|I_n|$ and $$ \sup_{|I_n|\le 1} |I_n|\log  \frac2{|I_n|}<\infty$$, 
		we get $$\avint_Q|v_N-\zeta(x_1)\added{a(I_n)}|\le C(n)<\infty$$
		uniformly in $Q$ and $N$. \added{The mean oscillation is at most twice this bound, and ball and cube BMO seminorms are equivalent.} 
	\end{proof}
	
	Using Lemma \ref{tests} and the above test function we obtain the following. 
	
	\begin{lemma} \label{affine}
		Given any $x\in E_f$, if $Df(x)\ne 0$, then $Df(x)$ is invertible and 
		$$ |Df(x)|\le  C\min_{\xi\in \added{S^{n-1}}} |Df(x)\xi|$$
		for some  constant  $C (n,C_f)$. 
	\end{lemma}
	\begin{proof} 
		Fix arbitrary $u\in C^0(\R^n)\cap L^\infty(\R^n)$.  We claim that 
		\begin{align}\label{claim-2}
			\avint_B |u\circ \ell_x -(u\circ \ell_x )_B| \le\added{C_f}[u]_{\BMO}\quad\forall B\subset\R^n 
		\end{align}
		and hence    $[u\circ\ell_x]_{\BMO}\le\added{C_f}[u]_{\BMO}$. 
		
		Given any ball $B$, since  Lemma~\ref{rescaling} gives 
		$u\circ f_{x,r}\to u\circ \ell_x $  in $L^1 (B)$,
		we obtain 
		$$\avint_{\added{B}}|u\circ\ell_x-(u\circ\ell_x)_B|=\lim_{r\added{\downarrow0}}\avint_{\added{B}}|u\circ f_{x,r}-(u\circ f_{x,r})_B|.$$
		Note that 
		\begin{align*}(u\circ f_{x,r})_B
			&=\added{\frac1{|B|}}\int_B u(\frac1r[f(x+rz)-f(x)]) \,dz\\
			& 
			=\frac1{|\added{x+rB}|}\int_{\added{x+rB}} u(\frac1r[f(y)-f(x)]) \,dy\\
			&=[u(\tfrac1r[f(\cdot)-f(x)])]_{x+rB}\added{=:c(B,x,r)}  
		\end{align*}
		Therefore  
		\begin{align*}
			\avint_B|u\circ f_{x,r}-(u\circ f_{x,r})_B|&= \avint_B|u\circ f_{x,r}-c(B,x,r)|\\
			&\added{=}\avint_{x+rB}|u(\tfrac1r[f(y)-f(x)])-[u(\tfrac1r[f(\added{\cdot})-f(x)])]_{x+rB}|\,\added{dy}   \\
			& \le C_f[u(\tfrac1r[\cdot-f(x)])]_{\BMO} 
		\end{align*}
		Since translations and scalar dilations preserve the BMO seminorm,
		we obtain 
		$$ [u(\tfrac1r[\cdot-f(x)])]_{\BMO} 
		=  [u]_{\BMO}. $$
		This gives the desired claim. 
		
		Up to  \reviewfix{orthogonal changes} of coordinates and a nonzero scalar normalization,
		we may  assume
		\[
		Df(x)=\operatorname{diag}(1,s_2,\ldots,s_n),
		\qquad1\ge s_2\ge\cdots\ge s_n\ge0.
		\]
		To be precise, \added{by singular value decomposition there are a scalar multiple $P$ of an orthogonal matrix and an orthogonal matrix $Q$} such that 
		$$PDf(x)Q=\operatorname{diag}(1,s_2,\ldots,s_n),
		\qquad1\ge s_2\ge\cdots\ge s_n\ge0.$$
		Let $\added{\widetilde f(z)=Pf(x+Qz)}$. We have 
		$\added{D\widetilde f(\added{0})=PDf(x)Q= \operatorname{diag}(1,s_2,\ldots,s_n)}$.
		Thus $$u\circ\widetilde f= u\circ P\circ f\circ\added{(x+Q\cdot)},$$
		and hence 
		$$[u\circ\widetilde f]_{\BMO}=[ u\circ P\circ f\circ (x+Q\cdot)]_{\BMO}=[ u\circ P\circ f]_{\BMO}\le C_f [u\circ P]_{\BMO}=C_f[u]_{\BMO}$$
		We may therefore replace $f,x$ by $\widetilde f,0$ below.

		If $s_n=0$, then on $Q=(-1,1)^n$ we have
		$v_N(Df(x)z)=Nz_1$. Thus
		$(v_N\circ\ell_x)_Q=0$
		and 
		$$  \avint_{Q}|v_N\circ\ell_x-(v_N\circ\ell_x)_Q|=\frac1{2 }\int_{-1}^1N|z_1|\,dz_1= \frac N2.$$
		By  
		Lemma \ref{tests} and \eqref{claim-2},  one has 
		$$\frac N2\le C(n) [v_N\circ\ell_x]_{\BMO}\le  C(n)C_f [v_N ]_{\BMO}\le C(n)C_f,$$
		which is \reviewfix{a contradiction} when $N$ is sufficiently large.

		If $s_n>0$,  on $Q=(-1,1)^n$, we have 
		\[
		v_N(Df(x)z)=z_1\min\{  \log\frac1{s_n|z_n|},N\}.
		\]
		Again the average $(v_N\circ \ell_x)_Q=0$. As $N\to\infty$, 
		\[\lim_{N\to\infty}\avint_Q|v_N\circ\ell_x| =\frac12\int_0^1\log\frac1{s_nt}\,dt=
		\frac12(1+ \log \frac1{s_n} ).
		\]
		The uniform bound from \added{\eqref{claim-2}} and Lemma \ref{tests} therefore bounds
		$\log \frac1{s_n}$ in terms of $n,C_f$.
	\end{proof}
	We are ready to prove Theorem \ref{main}.
	We use the standard fact that the Jacobian of a $W^{1,n}_{\loc}$
	homeomorphism has the sign of its topological orientation almost everywhere;
	see \cite[Theorems 1.1 and 1.2]{HenclMaly}.
	\begin{proof}[Proof of Theorem~\ref{main}]\label{mainproof}
		By Lemma~\ref{critical},
		$f\in W^{1,n}_{\loc}(\mathbb R^n)$.
		At $x\in E_f$, if $Df(x)=0$, the inequality $|Df(x)|^n\le  |\det Df(x)|$ is trivial.
		Otherwise, Lemma \ref{affine} tells $$ |Df(x)|\le K\min_{ \xi\in S^{n-1}}|Df(x)\xi|\le K|\det Df(x)|^{1/n},$$ 
		as desired, where $K$ depends only on $n,C_f$.
		Thus $|Df|^n\le K^n|\det Df|$ almost everywhere. After an
		orientation reversing reflection if necessary, $\det Df\ge0$ almost everywhere.
		The analytic definition of quasiconformality gives the conclusion.
	\end{proof}
	
	\section*{Declarations}
	
	\noindent{\bf Competing interests} 
	The authors have no competing interests to declare that are relevant to the content of this article.
	
		\noindent\textbf{Funding.} Yepei Wen and Yuan Zhou were supported by the
	National Natural Science Foundation of China (No.~12671115) and the
	Fundamental Research Funds for the Central Universities.
	
	\noindent{\bf Data availability}
	No datasets were generated or analysed during the current study.
	
\noindent{\bf Author address}
School of Mathematical Sciences, Beijing Normal University, Beijing 100875, P. R. China.

\noindent{\bf E-mail}
\texttt{202331130034@mail.bnu.edu.cn}, \texttt{yuan.zhou@bnu.edu.cn}

\end{document}